\documentclass[12pt]{article}

\usepackage{graphicx,color}
\usepackage{amsmath,amssymb,amsthm}
\renewcommand{\qed}{\hspace*{\fill}\blacksquare}
\newtheorem{thm}{Theorem}
\newtheorem{lem}[thm]{Lemma}
\newtheorem{cor}[thm]{Corollary}

\newtheorem{prop}[thm]{Proposition}

\newtheorem{defi}{Definition}

\newtheorem{notation}{Notation}

\def\beq{\begin{equation}}\def\eeq{\end{equation}}
\def\beqn{\begin{eqnarray}}\def\eeqn{\end{eqnarray}}

\def\qed{\ifhmode\unskip\nobreak\fi\quad\ifmmode\Box\else$\Box$\fi}

\title{The Fractional Haemers Bound of the Mycielski Construction\thanks{A short version of this work will appear as an extended abstract in the proceedings of Eurocomb 2025. The project supported by the Doctoral Excellence Fellowship Programme (DCEP) is funded
by the National Research Development and Innovation Fund of the Ministry of Culture and
Innovation and the Budapest University of Technology and Economics, under a grant
agreement with the National Research, Development and Innovation Office.}}
\author{\hfil Bence Csonka\thanks
  {Department of Computer Science and Information Theory,
Faculty of Electrical Engineering and Informatics,
Budapest University of Technology and Economics, E-mail: {\tt csonkab@edu.bme.hu}
 }}

\date{}

\begin{document}

\maketitle
\begin{abstract}
\par\noindent
We investigate the effect of the generalized Mycielski construction $M_r(G)$ on the complementary fractional Haemers bound $\bar{\mathcal{H}}_f(G; \mathbb{F})$, a parameter that depends on a graph $G$ and a field $\mathbb{F}$. The effect of the Mycielski construction on graph parameters has already been studied for the fractional chromatic number $\chi_f$ and the complementary Lov\'asz theta number $\bar{\vartheta}$. Larsen, Propp, and Ullman provided a formula for \( \chi_f(M_2(G)) \) in terms of $\chi_f(G)$. This was later generalized by Tardif to \( \chi_f(M_r(G)) \) for any $r$, and Simonyi and the author gave a similar expression for \( \bar{\vartheta}(M_2(G)) \) in terms of $\bar{\vartheta}(G)$. In this paper, we show that Tardif’s formula for the fractional chromatic number remains valid for \( \bar{\mathcal{H}}_f \) whenever $ \bar{\mathcal{H}}_f(G; \mathbb{F})$ equals the clique number of $G$. In particular, we provide a general upper bound on $\bar{\mathcal{H}}_f(M_r(G); \mathbb{F})$ in terms of $\bar{\mathcal{H}}_f(G;\mathbb{F})$ and we prove that this bound is tight whenever $ \bar{\mathcal{H}}_f(G; \mathbb{F})$ equals the clique number of $G$. Using the bounds, we present a general class of graphs for which the fractional Haemers bound of the generalized Mycielski construction can be determined exactly.
\end{abstract}

\section{Introduction}

\bigskip
\par\noindent
The study of graph parameters that provide upper bounds on Shannon capacity has been a central topic in zero-error information theory and extremal combinatorics since Shannon's foundational work~\cite{Sha56}.
\smallskip
\par\noindent
To introduce the Shannon OR-capacity, we define the OR-product of two graphs as follows. Let $G$ and $H$ be graphs, their OR-product $G\cdot H $ is defined by
\begin{gather*}
V(G \cdot H) = V(G) \times V(H); \\
E(G \cdot H) = \left\{ \{(g,h),(g',h')\} : \{g,g'\} \in E(G) \text{ or } \{h,h'\} \in E(H) \right\}.
\end{gather*}
\par\noindent
The Shannon OR-capacity is defined by $C_{\text{OR}}(G) = \lim_{t \to \infty} \sqrt[t]{\omega(G^t)}$,
where $G^t$ is the $t$-fold OR-product of $G$ with itself and $\omega(G)$ is the clique number of $G$. 
By Fekete's lemma, the limit exists. It is immediate from the definition that
$\omega(G)^t \leq \omega(G^t) \leq |V(G)|^t$,
which implies that the Shannon OR-capacity of any graph is always finite. For simplicity, in this paper we refer to the Shannon OR-capacity as the Shannon capacity.
\smallskip
\par\noindent
Shannon~\cite{Sha56} originally used a complementary formulation when defining graph capacity. The above definition is, however, more adequate for our discussion, just as is the case in Chapter 11 of the book~\cite{Imre11}.
In Shannon's original formulation, the definition uses the complementary language (AND-product), where two vertices are adjacent in the $t$-th AND-power if they are adjacent or equal in every coordinate and they differ in at least one coordinate. The corresponding AND-version of the capacity is denoted by $C_{\rm AND}(G)$. The two formulations are equivalent, since for the AND-version of the capacity we have the identity
$C_{\rm OR}(\bar{G}) = C_{\rm AND}(G)$, where $\bar{G}$ denotes the complement of $G$.
\medskip
\par\noindent
While the Shannon capacity is well-defined, it is difficult to compute. Lov\'asz~\cite{LL79} introduced the complementary theta number $\bar{\vartheta}(G)$, which provides a semidefinite programming-based upper bound on $C_{\rm OR}(G)$. The parameter $\bar{\vartheta}(G)$ inspired the study of graph parameters that behave similarly to the complementary theta number.
\medskip
\par\noindent
 To formalize this approach, Zuiddam~\cite{Zuiddam} introduced the asymptotic spectrum of graphs, denoted by $\mathcal{A}$. Before defining it, we recall the notion of graph homomorphism. Let $G$ and $H$ be two graphs. The function $f:V(G)\to V(H)$ is called a homomorphism if for all $v,w \in V(G)$, whenever $\{v,w\} \in E(G)$, we have $\{f(v),f(w)\} \in E(H)$. 
\smallskip
\par\noindent
Let \( \mathcal{G} \) denote the class of all finite graphs. The asymptotic spectrum $\mathcal{A}$ is the class of graph parameters \( \varphi: \mathcal{G} \to \mathbb{R}_{\geq 0} \) that satisfy the following properties:
\begin{enumerate}
    \item for any graphs $G$ and $H$, we have $\varphi(G+H) = \varphi(G) + \varphi(H)$, where $G+H$ denotes the graph obtained by taking disjoint copies of $G$ and $H$ and adding all possible edges between them;
    \item if there exists a graph homomorphism from $G$ to $H$, then $\varphi(G) \leq \varphi(H)$;
    \item for any graphs $G$ and $H$, we have $\varphi(G \cdot H) = \varphi(G) \varphi(H)$, where $\cdot$ denotes the OR-product of graphs;
    \item $\varphi(K_1) = 1$, where $K_1$ is the graph with one single vertex.
\end{enumerate}
Importantly, each graph parameter in $\mathcal{A}$ provides an upper bound on the Shannon capacity, that is, $C_{\rm OR}(G) \leq \varphi(G)$ for all $\varphi \in \mathcal{A}$.
\smallskip
\par\noindent
It is known that both $\chi_f$ and $\bar{\vartheta}$ belong to the asymptotic spectrum (i.e. $\chi_f, \bar{\vartheta} \in \mathcal{A}$).
\smallskip
\par\noindent
The following theorem of Zuiddam reveals the strong connection between Shannon capacity and the asymptotic spectrum. 
\smallskip
\par\noindent
For any graph $G$,
\[
C_{\rm OR}(G) = \min_{\varphi \in \mathcal{A}} \varphi(G).
\]
\par\noindent
That is, the Shannon capacity of a graph equals the minimum of the values assigned to it by the parameters in the asymptotic spectrum.
\medskip
\par\noindent
In recent decades, various authors have studied the connection between the asymptotic spectrum and the Mycielski construction, which is defined as follows.

\medskip
\par\noindent
\begin{defi}[\cite{Myc}]
  Let $G$ be a simple graph and $r$ a positive integer.
  Its generalized Mycielskian $M_r(G)$ is defined on the vertices 
  $$V(M_r(G))=V(G)\times\{0,1,...,r-1\}\cup \{z\}$$
  with edge set 
\begin{gather*}
    E(M_r(G)) = \{\{(v,0),(w,0)\}:\{v,w\} \in E(G)\} \cup \\
    \cup\{\{(v,i),(w,j)\}:\{v,w\} \in E(G) \text{ and } |i-j|=1\} \cup\{\{(v,r-1),z\}:v \in V(G)\}.
\end{gather*}
\end{defi}
\smallskip
\par\noindent
The base Mycielski construction $M_2$ demonstrates that the gap between the clique number and the chromatic number can be arbitrarily large. In particular, for every graph $G$, $\omega(M_2(G))=\omega(G)$ and $\chi(M_2(G)) = \chi(G)+1$, where $\chi(G)$ denotes the chromatic number of $G$. Additionally, any odd cycle can be obtained as $M_r(K_2) = C_{2r+1}$.
\smallskip
\par\noindent
Larsen, Propp, and Ullman were the first to establish a formula describing the relationship between the base Mycielski construction $M_2$ and a parameter in the asymptotic spectrum~\cite{LPU}. They showed that
\[
\chi_f(M_2(G)) = \chi_f(G) + \frac{1}{\chi_f(G)}.
\]
Later, this result was generalized by Tardif~\cite{CTcone} to any positive integer $r$.
\[
\chi_f(M_r(G)) = \chi_f(G) + \frac{1}{\sum_{k=0}^{r-1}(\chi_f(G)-1)^k}.
\] 
In~\cite{Csonka} it has been shown that a formula can also be provided for $\bar{\vartheta}(M_2(G))$ as a function of $\bar{\vartheta}(G)$, that is, the value of $\bar{\vartheta}(G)$ determines $\bar{\vartheta}(M_2(G))$. This result revealed that
\[
\bar{\vartheta}\left(M_2(G)\right) = \frac43\bar{\vartheta}(G)\cos{\left(\frac13{\rm arccos}\left(1 - \frac{27}{4\bar{\vartheta}(G)}+\frac{27}{4\bar{\vartheta}^2(G)}\right)\right)}-\frac13\bar{\vartheta}(G)+1.
\]
\par\noindent
\smallskip
The above theorems suggest the following natural question. If $\varphi \in \mathcal{A}$, then is $\varphi(M_r(G))$ determined by $\varphi(G)$? See Problem 5 in the open problems section of~\cite{Csonka}.
\medskip
\par\noindent
In the same paper where Lov\'asz introduced the complementary theta number, he asked whether $C_{\rm OR}(G) = \bar{\vartheta}(G)$ holds for every graph $G$. Haemers~\cite{Haemers,Haemers2} provided a negative answer to this question by introducing the complementary Haemers bound $\bar{\mathcal{H}}(G; \mathbb{F})$, which is defined for a graph $G$ and a field $\mathbb{F}$. In particular, he showed that for the complementary Haemers bound and the Schläfli graph $S$, it holds that
\[
\bar{\mathcal{H}}(S; \mathbb{R}) < \bar{\vartheta}(S).
\]
Blasiak \cite{Anna} introduced a generalization of the complementary Haemers bound, which was further investigated by Bukh and Cox~\cite{Boris}. 
\smallskip
\par\noindent
In the present paper, we focus on the complementary fractional Haemers bound $\bar{\mathcal{H}}_f(G;\mathbb{F})$. The following version of Schrijver's definition, as presented in~\cite{Boris}, is used in this paper.
\medskip
\par\noindent
\begin{defi}
    Let $\mathbb{F}$ be a field, and let $G$ be a graph. A collection $\{X_v\subseteq \mathbb{F}^n:v \in V(G)\}$ is called a dual $(n,d)$-representation of the graph $G$ over $\mathbb{F}$ if, for any vertex $v \in V(G)$, the following conditions hold:
    \begin{enumerate}
        \item $X_v$ is a subspace of $\mathbb{F}^n$; 
        \item $\dim{X_v}=d$; 
        \item $X_v \cap \left(\sum_{w \in N(v)}X_w\right) = \{0\}$, where the sum denotes the subspace generated by them, and $N(v)$ denotes the neighborhood of vertex $v$.
    \end{enumerate}
    The complementary fractional Haemers bound is defined as follows.
    \[
    \bar{\mathcal{H}}_f(G;\mathbb{F}) = \inf{\left\{\frac{n}{d}:G \text{ has a dual $(n,d)$-representation over $\mathbb{F}$}\right\}}.
    \]
\end{defi}
\smallskip
\par\noindent
The (non-fractional) complementary Haemers bound~\cite{Haemers} is defined by
\[\bar{\mathcal{H}}(G;\mathbb{F}) = \min\{n:G\text{ has a dual } (n,1)\text{-representation over }\mathbb{F} \}.
\]
Notice that $\bar{\mathcal{H}}_f(G;\mathbb{F}) \leq \bar{\mathcal{H}}(G;\mathbb{F})$ for every graph $G$ and for every field $\mathbb{F}$, since the definition of $\bar{\mathcal{H}}$ uses dual $(n,1)$-representation and the definition of $\bar{\mathcal{H}}_f$ uses 
a more general $(n,d)$-representation.
\smallskip
\par\noindent
Bukh and Cox showed that $\bar{\mathcal{H}}_f$ is an element of the asymptotic spectrum. Among the best-known elements of the spectrum are the fractional chromatic number $\chi_f$, the complementary theta number $\bar{\vartheta}$, and the fractional Haemers bound $\bar{\mathcal{H}}_f$. However, the (non-fractional) complementary Haemers bound is not in the spectrum, similarly to how the (non-fractional) chromatic number is not.
\smallskip
\par\noindent
Note that the fractional chromatic number and the complementary Lov\'asz theta number admit dual formulations; that is, both $\chi_f$ and $\bar{\vartheta}$ have equivalent minimum and maximum formulations, which are used in the proofs of the theorems about $\chi_f(M_r(G))$ and $\bar{\vartheta}(M_2(G))$. However, the existence of a dual formulation for $\bar{\mathcal{H}}_f$ is still unknown.
\smallskip
\par\noindent
Our next result shows that the formula given by Tardif for the fractional chromatic number provides a general upper bound for $\bar{\mathcal{H}}_f(M_r(G); \mathbb{F})$.
\smallskip
\par\noindent
\begin{thm}\label{Theorem1}
For any simple graph $G$ and any field $\mathbb{F}$, the following inequality holds:
\[
    \bar{\mathcal{H}}_f(M_r(G); \mathbb{F}) \leq \bar{\mathcal{H}}_f(G; \mathbb{F}) + \frac{1}{\sum_{k=0}^{r-1} (\bar{\mathcal{H}}_f(G; \mathbb{F}) - 1)^k}.
\]
\end{thm}
\smallskip
\par\noindent
According to Theorem~\ref{Theorem1}, for the Schläfli graph $S$, we have $\bar{\mathcal{H}}_f(M_r(S);\mathbb{R}) \leq 7 + \frac{6}{7^r-1}$ since Haemers proved that $\bar{\mathcal{H}}(S;\mathbb{R}) \leq 7$, we also have $\bar{\mathcal{H}}_f(S;\mathbb{R}) \leq 7$.
\smallskip
\par\noindent
This naturally raises the question of when the upper bound is tight. The following result shows that the bound is tight whenever $\bar{\mathcal{H}}_f(G; \mathbb{F}) = \omega(G)$.
\smallskip
\par\noindent
\begin{thm}\label{Theorem2}
For any field $\mathbb{F}$, positive integers $m,r$,
\[
    \bar{\mathcal{H}}_f(M_r(K_m);\mathbb{F}) = m + \frac{1}{\sum_{k=0}^{r-1}(m-1)^k},
\]
where $K_m$ is the complete graph with $m$ vertices.
\end{thm}
\smallskip
\par\noindent
\begin{cor}
    If $\bar{\mathcal{H}}_f(G;\mathbb{F}) = \omega(G)$, then the upper bound of Theorem~\ref{Theorem1} is tight.
\end{cor}
\smallskip
\par\noindent
 For example, it is true for perfect graphs.
 \smallskip
 \par\noindent
The following application shows that the bound can also be exact for graphs beyond the class of perfect graphs. Bukh and Cox~\cite{Boris} proved that for every field $\mathbb{F}$ of nonzero characteristic, there exists a graph $G$ such that for every field $\mathbb{F}'$,
\[
\bar{\mathcal{H}}_f(G;\mathbb{F}) < \min\{\bar{\mathcal{H}}(G;\mathbb{F}'),\bar{\vartheta}(G)\}.
\]
In their proof, they used the following graph family introduced by Haemers: let $J_n^p$ be the graph whose vertex set is $\binom{[n]}{p+1}$, where two sets $X,Y$ are adjacent if $|X \cap Y| \not\equiv 0 \pmod{p}$. If $(p+2)$ divides $n$ and $\mathbb{F}$ is a field of characteristic $p$, then $\bar{\mathcal{H}}(\bar{J_n^p};\mathbb{F}) = \omega(\bar{J_n^p}) = n$, however $\bar{\vartheta}(\bar{J_n^p}) = \left(\frac{p}{(p+1)^2}+o(1)\right)n^2$.
\smallskip
\par\noindent
Applying Theorem~\ref{Theorem1} and Theorem~\ref{Theorem2} to this graph, we obtain the following corollary.
\smallskip
\par\noindent
\begin{cor}
If $(p+2)$ divides $n$ and $\mathbb{F}$ is a field of characteristic $p$, then
\[
\bar{\mathcal{H}}_f\bigl(M_r(\bar{J_n^p});\mathbb{F}\bigr)
= n + \frac{1}{\sum_{k=0}^{r-1}(n-1)^k}.
\]
\end{cor}
\medskip
\par\noindent
Considering the above theorems, the following natural question arises. Is it true that there is a function $g:\mathbb{R} \rightarrow \mathbb{R}$ such that for every graph $G$ and every field $\mathbb{F}$, $\bar{\mathcal{H}}_f(M_r(G);\mathbb{F}) = g(\bar{\mathcal{H}}_f(G;\mathbb{F}))$? In particular, does Tardif’s formula for the fractional chromatic number hold more generally for the complementary fractional Haemers bound?
\smallskip
\par\noindent
Motivated by the behavior of $\bar{\mathcal{H}}_f$ under the Mycielski construction for cliques and the general upper bound, we conjecture that the identity
\[
    \bar{\mathcal{H}}_f(M_r(G); \mathbb{F}) = \bar{\mathcal{H}}_f(G; \mathbb{F}) + \frac{1}{\sum_{k=0}^{r-1} (\bar{\mathcal{H}}_f(G; \mathbb{F}) - 1)^k}
\]
holds for every graph $G$. We are not aware of any counterexample, even in cases where the complementary fractional Haemers bound is non-integer.
\medskip
\par\noindent
In Section 2, we prove Theorem~\ref{Theorem1} and in Section 3, we prove Theorem~\ref{Theorem2}.

\section{A General Upper Bound}

\bigskip
\par\noindent
To prove Theorem~\ref{Theorem1}, we use the concept of the tensor product of vector spaces.

\begin{defi}
Let \( v \in \mathbb{F}^n \) and \( w \in \mathbb{F}^m \). Their Kronecker product is defined as
\[
v \otimes w = (w_1 v, w_2 v, \dots, w_m v) \in \mathbb{F}^{nm},
\]
where each \( w_j v \in \mathbb{F}^n \) is a scalar multiple of \( v \).
    Let $V \subseteq \mathbb{F}^n$ and $W \subseteq \mathbb{F}^m$ be two subspaces. Let $\{v^{(1)},...,v^{(k)}\}$ be a basis of $V$ and let $\{w^{(1)},...,w^{(l)}\}$ be a basis of $W$. The tensor product of the two subspaces is the subspace $V \otimes W \subseteq \mathbb{F}^{nm}$ defined as
    \[
    V \otimes W = {\rm span}_{\mathbb{F}}\left\{v^{(i)}\otimes w^{(j)}:1 \leq i \leq k; 1 \leq j \leq l\right\},
    \]
\end{defi}
\medskip
\par\noindent
where ${\rm span}_\mathbb{F}$ denotes the generated subspace over $\mathbb{F}$. It is well-known that if $\{v^{(1)},...,v^{(k)}\} \subseteq \mathbb{F}^n$ is a linearly independent system and $\{w^{(1)},...,w^{(l)}\}\subseteq \mathbb{F}^m$ is linearly independent, then the collection of vectors $\{v^{(i)} \otimes w^{(j)}\}_{i,j}$ is also a linearly independent system.
\medskip
\par\noindent
The following identities can be found in Chapter 1 of~\cite{Ryan}. Although they are stated there for the fields $\mathbb{R}$ and $\mathbb{C}$, the same proofs remain valid over an arbitrary field $\mathbb{F}$. 
\medskip
\par\noindent
\begin{thm}\label{tensor}
    Let $\mathbb{F}$ be a field, and let $U,V,W \subseteq \mathbb{F}^n$ be subspaces. The following identities hold:
    \begin{enumerate}
        \item $\dim{(U \otimes V)} = \dim{U}\dim{V}$;
        \item $(U + V) \otimes W = U \otimes W + V \otimes W$;
        \item if $U \cap V = \{0\}$, then $(U \otimes W)\cap (V \otimes W) = \{0\}$.
    \end{enumerate}
\end{thm}
\medskip
\par\noindent
\smallskip
\par\noindent
Define the subspace $\Gamma [a,b] \subseteq \mathbb{F}^M$ for some positive integer $M$ and $1 \leq a \leq b \leq M$ as follows:
\begin{gather*}
    \Gamma [a,b] := {\rm span}_{\mathbb{F}}\{e^{(i)}:a \leq i \leq b\},
\end{gather*}
where the $e^{(i)}$'s are the standard basis vectors of the space $\mathbb{F}^M$.
\smallskip
\par\noindent
The following lemma presents a key structural property of subspaces $\Gamma[a,b]$.
\medskip
\par\noindent
\begin{lem}\label{intersection}
    Let $a\leq b$ be positive integers, and let $[a,b]:=\{s\in \mathbb{Z}_{\geq0}:a\leq s \leq b\}$ be called an integer interval. Let $M$ be a positive integer and $[a_1,b_1],\dots,[a_k,b_k],[c_1,d_1],\dots , [c_{\ell},d_{\ell}],[c,d] \subseteq [1,M]$ be integer intervals such that $\left(\bigcup_{i=1}^k [a_i,b_i]\right)\cap[c,d] = \emptyset$. If $X_1,...,X_k,Y_1,...,Y_{\ell},Y$ are subspaces of a vector space over $\mathbb{F}$, then
    \begin{gather*}
    \left[\sum_{i=1}^k\left(X_i\otimes \Gamma[a_i,b_i]\right)\right] \cap \left[\sum_{j=1}^{\ell}(Y_j\otimes\Gamma[c_j,d_j]) + Y \otimes \Gamma[c,d]\right] = \\
    = \left[\sum_{i=1}^k\left(X_i\otimes \Gamma[a_i,b_i]\right)\right] \cap \left[\sum_{j=1}^{\ell}(Y_j\otimes\Gamma[c_j,d_j]) \right].
    \end{gather*}
\end{lem}
\smallskip
\par\noindent
The proof of Lemma~\ref{intersection} can be found in Appendix A.
\medskip
\par\noindent
The following proof is inspired by Tardif’s construction for the fractional chromatic number~\cite{CTcone}.
\medskip
\par\noindent
{\em Proof of Theorem~\ref{Theorem1}.} Let $\{X_v \subseteq \mathbb{F}^n : v \in V(G)\}$ be a dual $(n,d)$-representation of $G$ over the field $\mathbb{F}$. If $r=1$, then $\bar{\mathcal{H}}_f(M_1(G);\mathbb{F}) = \bar{\mathcal{H}}_f(G + K_1;\mathbb{F}) = \bar{\mathcal{H}}_f(G;\mathbb{F}) + 1$, where we use $\bar{\mathcal{H}}_f \in \mathcal{A}$. The other trivial case $\bar{\mathcal{H}}_f(G;\mathbb{F})=1$ implies $G$ consists of isolated vertices, therefore $M_r(G)$ is a disjoint union of a star graph and isolated vertices, which implies $\bar{\mathcal{H}}_f(M_r(G);\mathbb{F}) = 2$. Our goal is to construct a dual $(N,D)$-representation of $M_r(G)$, where $r \geq 2$, $\bar{\mathcal{H}}_f(G;\mathbb{F}) \geq 2$ and
\begin{gather*}
    N := n\sum_{i=0}^{r-1}d^{2r-2-i}(n-d)^i + d^{2r-1};\\
    D := \sum_{i=0}^{r-1}d^{2r-1-i}(n-d)^i.
\end{gather*}
\smallskip
\par\noindent
A direct calculation shows that
\[
\frac{N}{D} = \frac{n}{d} + \frac{1}{\sum_{i=0}^{r-1}\left(\frac{n}{d}-1\right)^i},
\]
which is the statement of Theorem~\ref{Theorem1}. 
\medskip
\par\noindent
We construct an $(N,D)$-representation explicitly.
\smallskip
\par\noindent
For simplicity, we introduce the following notations:
\[
a_i := \sum_{j=0}^{i}d^{2r-2-j}(n-d)^j, \quad a_{-1}:= 0.
\]
\smallskip
\par\noindent
Let $M := \sum_{i=0}^{r-1}d^{2r-2-i}(n-d)^i + d^{2r-1}$ and all $\Gamma[a,b]$ be a subspace of $\mathbb{F}^M$.
\smallskip
\par\noindent
We assign to each vertex of $M_r(G)$ a subspace of $\mathbb{F}^n \otimes \mathbb{F}^M$. Since $\mathbb{F}^n \otimes \mathbb{F}^M \cong \mathbb{F}^{nM}$, each such subspace can be identified with a subspace of $\mathbb{F}^{nM}$.
\bigskip
\par\noindent
We define an assignment for each vertex $(v,k) \in V(M_r(G))$.
\medskip
\par\noindent
\[
    \tilde{X}_{(v,0)} := X_v \otimes \Gamma [1,a_{r-1}].
\]
\smallskip
\par\noindent
Let $E := {\rm span}_{\mathbb{F}}\{e^{(1)}\}$, where $e^{(1)} \in \mathbb{F}^n$ is the canonical vector whose first coordinate is $1$ and all other coordinates are zero.
\[
    \tilde{X}_{(v,1)} :=
    X_v \otimes\Gamma[d^{2r-2}+1,a_{r-1}] + E\otimes\Gamma\left[a_{r-1}+1,a_{r-1}+d^{2r-1}\right].
\]
If $k = 2\ell+1$,

\smallskip
\par\noindent
\footnotesize
\begin{gather*}
    \tilde{X}_{(v,k)} :=\\
    =\left(\sum_{w \in V(G)}X_w\right) \otimes\left(\sum_{i=0}^{\ell-1} \Gamma[a_{2i}+1,a_{2i+1}]\right) + X_v \otimes \Gamma[a_{k-1}+1,a_{r-1}] + E\otimes \Gamma\left[a_{r-1}+1,a_{r-1}+d^{2r-1}\right].
\end{gather*}
\normalsize
\smallskip
\par\noindent
If $k = 2\ell$,
\smallskip
\par\noindent
\begin{gather*}
    \tilde{X}_{(v,k)} := \left(\sum_{w \in V(G)}X_w\right) \otimes \left(\sum_{i=0}^{\ell-1} \Gamma[a_{2i-1}+1,a_{2i}]\right)+X_v \otimes \Gamma[a_{k-1}+1,a_{r-1}].
\end{gather*}

\smallskip
\par\noindent
We now define the assignment for the vertex $z$. The assignment is divided into two cases. 
\smallskip
\par\noindent
If $r = 2t$, then
\begin{gather*}
    \tilde{X}_z := \left(\sum_{v \in V(G)}X_v\right)\otimes\left(\sum_{i=0}^{t-1} \Gamma[a_{2i-1}+1,a_{2i}]\right).
\end{gather*}
If $r=2t+1$, then
\begin{gather*}
    \tilde{X}_z := \left(\sum_{v \in V(G)}X_v\right)\otimes\left(\sum_{i=0}^{t-1} \Gamma[a_{2i}+1,a_{2i+1}]\right) + E\otimes\Gamma\left[a_{r-1}+1,a_{r-1}+d^{2r-1}\right].
\end{gather*}
\smallskip
\par\noindent
We show that $\{\tilde{X}_{u}:u \in V(M_r(G))\}$ is a dual $(N,D)$-representation over $\mathbb{F}$.
\smallskip
\par\noindent
The next step is to verify in certain cases that the assignment indeed yields a valid dual representation.
\smallskip
\noindent
In the proof, we use the following notations. For all $v \in V(G)$ let $N_G(v)$ be the neighborhood of $v$ in $G$ and for all $v \in V(M_r(G))$ let $N_{M_r(G)}(v)$ be the neighborhood of $v$ in $M_r(G)$.
\smallskip
\par\noindent
Let $(v,0)$ be a vertex at the level $0$. Then
\smallskip
\par\noindent
\begin{gather*}
\tilde{X}_{(v,0)} \cap \left(\sum_{u \in N_{M_r(G)}((v,0))}\tilde{X}_{u}\right) = \\
= (X_v \otimes \Gamma[1,a_{r-1}]) \cap \left(\sum_{q \in N_G(v)}\tilde{X}_{(q,0)} + \sum_{q \in N_G(v)}\tilde{X}_{(q,1)}\right) \subseteq\\
\subseteq  (X_v \otimes \Gamma[1,a_{r-1}])
\cap \left[\left(\sum_{q \in N_G(v)}X_q\right)\otimes \Gamma[1,a_{r-1}] + E\otimes\Gamma\left[a_{r-1}+1,a_{r-1}+d^{2r-1}\right]\right].
\end{gather*}
Since $[1,a_{r-1}]\cap \left[a_{r-1}+1,a_{r-1}+d^{2r-1}\right] = \emptyset$, we can apply Lemma~\ref{intersection}, therefore it holds that
\[
\left(X_v \otimes \Gamma[1,a_{r-1}]\right) \cap \left[\left(\sum_{q \in N_G(v)}X_q\right) \otimes \Gamma[1,a_{r-1}]\right] = \{0\},
\]
using the fact $X_v \cap \left(\sum_{q \in N_G(v)}X_q\right) = \{0\}$ and Theorem~\ref{tensor}.
\smallskip
\par\noindent
The remaining cases can be found in Appendix B1.
\smallskip
\par\noindent
It remains to verify that all constructed subspaces have the correct dimension. The dimension of the whole space,
\begin{gather*}
    \dim{\left[\sum_{u \in V(M_r(G))}\tilde{X}_{u}\right]} = \dim{\left[\sum_{k=0}^{r-1}\sum_{v \in V(G)}\tilde{X}_{(v,k)} + \tilde{X}_z\right]} \leq \\
    \leq\dim{\left[\left(\sum_{v \in V(G)}X_v\right)\otimes \Gamma[1,a_{r-1}] + E \otimes\Gamma\left[a_{r-1}+1,a_{r-1}+d^{2r-1}\right]\right]} = \\ =n\sum_{i=0}^{r-1}d^{2r-2-i}(n-d)^i + d^{2r-1} = N. 
\end{gather*}
Finally, we verify that \(\dim \tilde{X}_{(v, k)} = D\) holds for all \((v, k) \in V(M_r(G))\) and \(\dim \tilde{X}_z = D\). These calculations can be found in Appendix B2.
\medskip
\par\noindent
We have shown that $M_r(G)$ has an $(N,D)$-representation over $\mathbb{F}$, hence $\bar{\mathcal{H}}_f(M_r(G);\mathbb{F}) \leq \frac{N}{D}$ follows.
\hfill $\qedsymbol$
\section{Tightness for $K_m$}

\par\noindent
In this section, we establish Theorem~\ref{Theorem2} by constructing a lower bound for $\bar{\mathcal{H}}_f(M_r(K_m); \mathbb{F})$. 
\smallskip
\par\noindent
It should be noted that the inclusion-exclusion principle is not true for $\dim{\left(\sum_{i}X_i\right)}$ in general, since the intersection of subspaces is not distributive over the addition of subspaces. For example, if $X,Y,Z\subseteq \mathbb{R}^2$ are different lines, then $(X+Y)\cap Z = Z$ since $X+Y = \mathbb{R}^2$, but $X\cap Z + Y \cap Z = \{0\}$. Fortunately, the inclusion-exclusion principle is valid for two subspaces. The following theorem is the so-called Grassmann's identity.
\smallskip
\par\noindent
Theorem~\ref{grassman} can be found in Chapter 2 Theorem 6 of the book~\cite{Hoffman}.
\smallskip
\par\noindent
\begin{thm}\label{grassman}
    For every subspace $X,Y \subseteq\mathbb{F}^n$ \[
    \dim{(X + Y)} = \dim{X} + \dim{Y} - \dim{(X \cap Y)}.
    \]
\end{thm}
\medskip
\par\noindent
\begin{notation}
    Let $A \subseteq V(G)$ be an arbitrary set of vertices and let $X_A \subseteq \mathbb{F}^n$ be the following subspace:
    \[
    X_A = \sum_{v \in A}X_v.
    \]
    If $A,B \subseteq V(G)$ are two disjoint vertex sets such that every vertex from $A$ is adjacent to every vertex from $B$, we denote this by $A \sim B$.
\end{notation}
\medskip
\par\noindent
\begin{prop}\label{independent}
    Let $A,B \subseteq V(G)$ be vertex sets such that $A \sim B$ and $B$ forms a clique, then $\dim{(X_A + X_B)} = \dim{X_A} + |B|d$. 
\end{prop}
\smallskip
\par\noindent
{\em Proof of Proposition~\ref{independent}.} We prove the proposition by induction on the size of $|B|$. First, consider the case where $B:=\{b\}$, that is, $|B| = 1$. Since the sum of two subspaces satisfies the inclusion-exclusion formula according to Theorem~\ref{grassman}, we can use it to decompose the dimensions.
\[
\dim{(X_A + X_b)} = \dim{X_A} + d - \dim{(X_A \cap X_b)}.
\]
According to the assumption, $A \sim B = \{b\}$, which implies that $X_A \subseteq \sum_{v \in N(b)} X_v$, from which we obtain $X_A \cap X_b = \{0\}$.

\smallskip
\par\noindent
Now, assume that the statement holds for $|B| = k$ and consider the case $|B| = k+1$ such that $B$ is a clique and $A \sim B$. We write $X_B$ as
\[
X_B = X_{B\setminus \{b\}} + X_b.
\]
Applying Theorem~\ref{grassman},
\[
\dim{\left(X_A + X_{B \setminus{\{b\}}} + X_b\right)} = \dim{\left(X_A + X_{B \setminus{\{b\}}}\right)} + \dim{X_b} - \dim{\left(\left(X_A + X_{B \setminus{\{b\}}}\right) \cap X_b\right)}.
\]
Since $B$ forms a clique and $X_A \subseteq X_{N(b)}$, it follows that $\left(X_A + X_{B \setminus{\{b\}}}\right) \cap X_b = \{0\}$. Using the induction hypothesis for $\dim{\left(X_A + X_{B \setminus{\{b\}}}\right)}$, we obtain
\[
\dim{\left(X_A + X_{B \setminus{\{b\}}} + X_b\right)} = \dim{X_A} + kd+d=\dim{X_A} +(k+1)d.
\]

\hfill $\qedsymbol$
\smallskip
\par\noindent
The proof of Theorem~\ref{Theorem2} is an extension of the proof for odd cycles in \cite{Boris}.
\bigskip
\par\noindent
{\em Proof of Theorem~\ref{Theorem2}.} As we mentioned in the beginning of the proof of Theorem~\ref{Theorem1}, the interesting case is $r\geq 2$ and $m \geq 2$. Let $\mathbb{F}$ be a field and $m \geq 2$. We know that $\bar{\mathcal{H}}_f(K_m;\mathbb{F}) = \chi_f(K_m)=m$, so we can apply Tardif's theorem~\cite{CTcone}, which implies that for any $r \geq 2$
\[
\chi_f(M_r(K_m)) = m + \frac{1}{\sum_{k=0}^{r-1}(m-1)^k}.
\]
According to~\cite{Boris}, $\bar{\mathcal{H}}_f(M_r(K_m);\mathbb{F}) \leq \chi_f(M_r(K_m))$, therefore the upper bound follows.

\smallskip
\par\noindent
We first prove the lower bound for $r=2$ and $r=3$. These cases illustrate the main ideas of the argument. We then turn to general $r$, where the proof is based on recursive estimates for suitable dimensions.

\smallskip
\par\noindent
\textbf{\underline{$r=2$:}} We prove that
\[
m + \frac{1}{m} \leq \bar{\mathcal{H}}_f(M_2(K_m);\mathbb{F}).
\]
Let $\{X_{(i,k)}:1\leq i \leq m;k=0,1\}\cup\{X_z\}$ be an arbitrary dual $(n,d)$-representation of the graph $M_2(K_m)$ over $\mathbb{F}$. For any $i=1,..,m$, the following is true:
\begin{equation}
    \dim{\left[X_{(i,0)}\cap X_{(i,1)}\right]} \geq (m+1)d-n.\label{eq5}
\end{equation}
To prove this inequality, consider the following inequality:
\begin{equation}
   n \geq \dim{\left[\sum_{j=2}^m X_{(j,0)}+X_{(1,0)}+X_{(1,1)}\right]}.\label{eq0} 
\end{equation}
The set of vertices $\{(2,0),...,(m,0)\}$ forms a clique and each vertex is adjacent to the vertices $(1,0),(1,1)$, and thus according to Proposition~\ref{independent},
\[
\dim{\left[\sum_{j=2}^m X_{(j,0)}+X_{(1,0)}+X_{(1,1)}\right]} = (m-1)d + \dim{\left[X_{(1,0)}+ X_{(1,1)}\right]}.
\]
Let us apply Theorem~\ref{grassman},
\[
\dim{\left[X_{(1,0)}+ X_{(1,1)}\right]} = \dim{X_{(1,0)}} + \dim{X_{(1,1)}} - \dim{\left[X_{(1,0)}\cap X_{(1,1)}\right]}.
\]
Substituting this into the inequality~\eqref{eq0}, we obtain the desired bound. Throughout the proof, we repeatedly use the symmetry among the vertices $(i,0)$, which allows us to assume $i=1$ without loss of generality.

\smallskip
\par\noindent
Now, we consider the following inequality:
\begin{equation}
n \geq \dim{\left[\sum_{j=2}^m \left(X_{(j,0)}\cap X_{(j,1)}\right)+X_{(1,1)} + X_z\right]}.\label{eq1}
\end{equation}
The vertex $z$ is adjacent to all $(1,1),...,(m,1)$ and the vertex $(1,1)$ is adjacent to all $(2,0),...,(m,0)$, hence firstly $X_z$ can be factored out of the sum, secondly $X_{(1,1)}$ can be factored out of the sum.
\[
\dim{\left[\sum_{j=2}^m \left(X_{(j,0)}\cap X_{(j,1)}\right)+X_{(1,1)} + X_z\right]} = \dim{\left[\sum_{j=2}^m \left(X_{(j,0)}\cap X_{(j,1)}\right) \right]} + 2d.
\]
Note that for any $j=2,\dots,m$, we have
\[
\left(X_{(j,0)}\cap X_{(j,1)}\right) \cap \left[\sum_{\substack{s=2 \\ s\neq j}}^m\left(X_{(s,0)}\cap X_{(s,1)}\right)\right] \subseteq X_{(j,0)}\cap \left[\sum_{\substack{s=2\\s\neq j}}^mX_{(s,0)}\right] \subseteq X_{(j,0)}\cap\left[\sum_{v \in N((j,0))}X_v\right] = \{0\}.
\]
Thus, for any $j$, we also have the following: 
\[
\left(X_{(j,0)}\cap X_{(j,1)}\right) \cap\left[\sum_{\substack{s=2 \\ s\neq j}}^m\left(X_{(s,0)}\cap X_{(s,1)}\right)\right] = \{0\},
\]
allowing us to extract the $X_{(j,0)}\cap X_{(j,1)}$ terms from the sum~\eqref{eq1}:
\begin{gather*}
n \geq \dim{\left[\sum_{j=2}^m \left(X_{(j,0)}\cap X_{(j,1)}\right)\right]} + 2d = \sum_{j=2}^m\dim{\left[X_{(j,0)}\cap X_{(j,1)}\right]} + 2d \geq\\
\geq (m-1)((m+1)d-n) + 2d.
\end{gather*}
Rearranging the inequality, we have
\[
\frac{n}{d}\geq m + \frac{1}{m}.
\]
\medskip
\par\noindent
\textbf{\underline{$r=3$:}} Now let $\{X_{(i,k)}:1\leq i \leq m;k=0,1,2\}\cup\{X_z\}$ be an arbitrary dual $(n,d)$-representation of the graph $M_3(K_m)$ over $\mathbb{F}$. Inequality~\eqref{eq5} still holds for any $i$:
\[
\dim{\left[X_{(i,0)} \cap X_{(i,1)}\right]} \geq (m+1)d -n.
\]
Using the bound above, we now show that for any $i$,
\begin{equation}
\dim{\left[X_{(i,0)}\cap X_{(i,2)}\right]}\geq (m^2-m+1)d-(m-1)n.\label{eq4}
\end{equation}
We write the following inequality:
\begin{equation}
n \geq \dim{\left[\sum_{j=3}^m \left(X_{(j,0)}\cap X_{(j,1)}\right) + X_{(1,0)} + X_{(1,2)} + X_{(2,1)}\right]}.\label{eq6}
\end{equation}
It is evident that
\[
X_{(2,1)} \cap \left[\sum_{j=3}^m \left(X_{(j,0)}\cap X_{(j,1)}\right) + X_{(1,0)} + X_{(1,2)}\right] = \{0\}.
\]
Therefore, $X_{(2,1)}$ can be factored out of the sum, leaving only the term
\[
\dim{\left[\sum_{j=3}^m \left(X_{(j,0)}\cap X_{(j,1)}\right) + X_{(1,0)} + X_{(1,2)} \right]}
\]
to be further decomposed. For arbitrary $j=3,..,m$, we have
\[
\sum_{\substack{s=2\\ s \neq j}}^m\left(X_{(s,0)}\cap X_{(s,1)}\right) + X_{(1,0)} + X_{(1,2)} \subseteq \sum_{v \in N((j,1))}X_v,
\]
so the terms $\sum_{j=3}^m \left(X_{(j,0)}\cap X_{(j,1)}\right)$ can be factored out of the sum~\eqref{eq6}.
\begin{gather*}
    n \geq \sum_{j=3}^m \dim{\left[X_{(j,0)}\cap X_{(j,1)}\right]} + \dim{\left[X_{(1,0)}+X_{(1,2)}\right]} + d \geq \\
    \geq (m-2)((m+1)d-n) + 3d - \dim{\left[X_{(1,0)}\cap X_{(1,2)}\right]}.
\end{gather*}
Hence the inequality~\eqref{eq4} follows.
\smallskip
\par\noindent
Finally, we now proceed analogously to the $r=2$ case, incorporating the subspace $X_z$:
\begin{gather*}
    n \geq \dim{\left[\sum_{j=1}^m \left(X_{(j,0)}\cap X_{(j,2)}\right) + X_z\right]}.
\end{gather*}
Since
\[
\sum_{j=1}^m \left(X_{(j,0)}\cap X_{(j,2)}\right) \subseteq \sum_{v \in N(z)}X_v,
\]
The term $X_z$ can be removed from the dimension count, since it has trivial intersection with the sum of the remaining subspaces. After that, the expression $\dim\left[\sum_{j=1}^m \bigl(X_{(j,0)} \cap X_{(j,2)}\bigr)\right]$
can be decomposed into a sum of dimensions, since every $(j,0)$ and $(j',0)$ are adjacent to each other.
\begin{gather*}
    n \geq \sum_{j=1}^m \dim{\left[X_{(j,0)}\cap X_{(j,2)}\right]} + d \geq m((m^2-m+1)d-(m-1)n) + d.
\end{gather*}
Rearranging for $\frac{n}{d}$, we obtain
\[
\frac{n}{d} \geq m + \frac{1}{m^2-m+1},
\]
which is the desired bound.
\medskip
\par\noindent
In the case $r \geq 4$, all terms needed for the dimension count appear, allowing the general recursion to be written down.
\medskip
\par\noindent
\textbf{\underline{$r \geq 4$:}} We provide a lower bound for $\bar{\mathcal{H}}_f(M_r(K_m);\mathbb{F})$, when $r \geq 4$. Similarly to the previous cases, consider an $(n,d)$-representation $\left\{X_{(i,k)}:1 \leq i \leq m;k=0,...,r-1\}\cup \{X_z\right\}$ over $\mathbb{F}$. We proceed by formulating bounds for each $2\leq k <r-1$ and demonstrating how new bounds can be derived for $k+1$ using the previous bounds. We distinguish two cases based on the parity of $k$: even and odd.

\medskip
\par\noindent
$\underline{\textbf{$k$ even}}:$ Let $k=2\ell$ and 
assume that the following inequalities hold for all $i=1,...,m$:
\begin{gather*}
    a_{\ell} \leq \dim{\left[\bigcap_{t=1}^{\ell}X_{(i,2t-1)}\right]}, \hspace{0.2 cm}
    b_{\ell} \leq \dim{\left[\bigcap_{t=0}^{\ell}X_{(i,2t)}\right]}, \hspace{0.2 cm}
    c_{2\ell-1} \leq \dim{\left[\bigcap_{t=0}^{2\ell-1}X_{(i,t)}\right]}.
\end{gather*}
We already have $a_0,a_1,b_0,c_0$, and $c_1$ from the proofs of the cases $r=2$ and $r=3$. Using the previous lower bounds, we derive a bound for $2\ell+1$, determining the estimates for $a_{\ell+1}$ and $c_{2\ell}$. (The value $b_{\ell+1}$ is not achieved in the even case.) 

\smallskip
\par\noindent
Firstly, we determine a suitable value for $c_{2\ell}$.

\smallskip
\par\noindent
\begin{gather*}
    n \geq \dim{\left[\sum_{j=2}^m X_{(j,0)} + \bigcap_{t=0}^{\ell}X_{(1,2t)} + \bigcap_{t=1}^{\ell} X_{(1,2t-1)}\right]} =\\
    =(m-1)d + \dim{\left[\bigcap_{t=0}^{\ell}X_{(1,2t)}\right]} + \dim{\left[\bigcap_{t=1}^{\ell} X_{(1,2t-1)}\right]} - \dim{\left[\bigcap_{t=0}^{2\ell}X_{(1,t)}\right]} \geq\\
    \geq (m-1)d + b_{\ell} + a_{\ell} -\dim{\left[\bigcap_{t=0}^{2\ell}X_{(1,t)}\right]}.
\end{gather*}
Thus, setting
\[
c_{2\ell} := (m-1)d +  b_{\ell} + a_{\ell} -n
\]
is an appropriate choice for the lower bound of the intersection dimension.

\smallskip
\par\noindent
As a second step, we derive a lower bound for $a_{\ell+1}$ using the previous inequalities. To this end, we consider the following inequality:
\begin{equation}
    n \geq \dim{\left[\sum_{j=3}^m\left(\bigcap_{t=0}^{2\ell}X_{(j,t)}\right) + X_{(1,2\ell+1)} + \bigcap_{t=1}^{\ell}X_{(1,2t-1)} + X_{(2,2\ell)}\right]}.\label{eq2}
\end{equation}
The subspace $X_{(2,2\ell)}$ has trivial intersection with the subspace generated by the other subspaces, so it can be factored out from the dimension. Furthermore, for any $j=3,..,m$, it holds that
\begin{gather*}
     \left[\bigcap_{t=0}^{2\ell} X_{(j,t)}\right]\cap \left[\sum_{\substack{s = 3\\ s \neq j}}\left(\bigcap_{t=0}^{2\ell}X_{(s,t)}\right) + X_{(1,2\ell+1)} + \bigcap_{t=1}^{\ell}X_{(1,2t-1)}\right] \subseteq \\
    \subseteq X_{(j,2\ell)} \cap \left[\sum_{\substack{s = 3\\ s \neq j}}X_{(s,2\ell-1)} + X_{(1,2\ell+1)} + X_{(1,2\ell-1)}\right] \subseteq X_{(j,2\ell)}\cap \left[\sum_{v \in N((j,2\ell))}X_v\right] = \{0\}.
\end{gather*}
Thus, inequality~\eqref{eq2} can be expanded as follows:
\begin{gather*}
    n \geq \sum_{j=3}^m\dim{\left[\bigcap_{t=0}^{2\ell}X_{(j,t)}\right]} + d + \dim{\left[\bigcap_{t=1}^{\ell}X_{(1,2t-1)}\right]} + d - \dim{\left[\bigcap_{t=1}^{\ell+1}X_{(1,2t-1)} \right]} \geq \\
    \geq (m-2)c_{2\ell} + a_{\ell} + 2d - \dim{\left[\bigcap_{t=1}^{\ell+1}X_{(1,2t-1)} \right]} .
\end{gather*}
From this we obtain that $a_{\ell+1}:=(m-2)c_{2\ell} + a_{\ell} + 2d-n$ serves as a valid lower bound.
\medskip
\par\noindent
$\underline{\textbf{$k$ odd}}:$ Let $k=2\ell+1$ and assume that we know the following inequalities for all $i=1,...,m$:
\begin{gather*}
    a_{\ell+1} \leq \dim{\left[\bigcap_{t=1}^{\ell+1}X_{(i,2t-1)}\right]}, \hspace{ 0.2 cm}
    b_{\ell} \leq \dim{\left[\bigcap_{t=0}^{\ell}X_{(i,2t)}\right]}, \hspace{0.2 cm}
    c_{2\ell} \leq \dim{\left[\bigcap_{t=0}^{2\ell}X_{(i,t)}\right]}.
\end{gather*}
We present an estimate for $c_{2\ell+1}$ as well as an estimate for $b_{\ell+1}$. (The value $a_{\ell+1}$ is not achieved in the odd case.) 

\smallskip
\par\noindent
Similarly to the previous even case, we write the following inequality:
\begin{gather*}
    n \geq \dim{\left[\sum_{j=2}^m X_{(j,0)} + \bigcap_{t=0}^{\ell}X_{(1,2t)} + \bigcap_{t=1}^{\ell+1} X_{(1,2t-1)}\right]} =\\
    =(m-1)d + \dim{\left[\bigcap_{t=0}^{\ell}X_{(1,2t)}\right]} + \dim{\left[\bigcap_{t=1}^{\ell+1} X_{(1,2t-1)}\right]} - \dim{\left[\bigcap_{t=0}^{2\ell+1}X_{(1,t)}\right]} \geq\\
    \geq (m-1)d + b_{\ell} + a_{\ell+1} -\dim{\left[\bigcap_{t=0}^{2\ell+1}X_{(1,t)}\right]}.
\end{gather*}
We thus define $c_{2\ell+1} = (m-1)d + a_{\ell+1} + b_{\ell} -n$, which serves as a valid lower bound. Next, we present an estimate for $b_{\ell+1}$.
\begin{gather*}
    n \geq \dim{\left[\sum_{j=3}^m\left(\bigcap_{t=0}^{2\ell+1}X_{(j,t)}\right) + X_{(1,2\ell+2)} + \bigcap_{t=0}^{\ell}X_{(1,2t)} + X_{(2,2\ell+1)}\right]}.
\end{gather*}
By analogous reasoning as in the previous case, the dimension of the sum can be expanded as follows.
\begin{gather*}
    n \geq \sum_{j=3}^m\dim{\left[\bigcap_{t=0}^{2\ell+1}X_{(j,t)}\right]} + d + \dim{\left[\bigcap_{t=0}^{\ell}X_{(1,2t)}\right]} - \dim{\left[\bigcap_{t=1}^{\ell+1}X_{(1,2t)} \right]} + d \geq \\
    \geq (m-2)c_{2\ell+1} + b_{\ell} + 2d - \dim{\left[\bigcap_{t=1}^{\ell+1}X_{(1,2t)} \right]} .
\end{gather*}
It follows that $b_{\ell+1}: = (m-2)c_{2\ell+1}  + b_{\ell} +2d-n$ serves as a valid lower bound.
\smallskip
\par\noindent
We now summarize the recursion relations derived above:
\begin{gather*}
    c_{2\ell} = (m-1)d + a_{\ell} + b_{\ell} -n;\\
    c_{2\ell+1} = (m-1)d + a_{\ell+1} + b_{\ell} -n; \\ 
    a_{\ell+1} = (m-2)c_{2\ell} + a_{\ell} + 2d-n; \\
    b_{\ell+1} = (m-2)c_{2\ell+1} + b_{\ell} + 2d-n.
\end{gather*}
Using the obtained recursions, we derive a formula for $c_k$.
\medskip
\par\noindent
\begin{lem}\label{lem2} For $k\geq 1$,
    \begin{gather*}
    c_{k} = d\left(m\sum_{i=0}^{k-1}(m-1)^i + 1\right) -n\sum_{i=0}^{k-1}(m-1)^i.
\end{gather*}
\end{lem}

\medskip
\par\noindent
 {\em Proof of Lemma~\ref{lem2}.} First, we show that the following recurrence holds for every $k \geq 1$:
\[
c_k = (m-1)c_{k-1} + 2d-n.
\]
If $k=2\ell$, then
\[
c_{2\ell}-c_{2\ell-1} = b_{\ell}-b_{\ell-1} = (m-2)c_{2\ell-1} + 2d-n,
\]
from which it follows that
\[
c_{2\ell} = (m-1)c_{2\ell-1} + 2d-n.
\]
Similarly, if $k=2\ell+1$, then
\[
c_{2\ell+1}-c_{2\ell} = a_{\ell+1}-a_{\ell} = (m-2)c_{2\ell} + 2d-n.
\]
After rearranging, we obtain the recurrence. Expanding the recurrence yields the following closed formula:
\[
c_k = (m-1)^k d + (2d-n)\sum_{i=0}^{k-1}(m-1)^i.
\]
Let us now simplify this closed formula.
\begin{gather*}
(m-1)^kd + (2d-n)\sum_{i=0}^{k-1}(m-1)^i = d\left(\sum_{i=0}^k(m-1)^i + \sum_{i=0}^{k-1}(m-1)^i\right) -n\sum_{i=0}^{k-1}(m-1)^i = \\
d\left((m-1)\sum_{i=0}^{k-1}(m-1)^i + 1 + \sum_{i=0}^{k-1}(m-1)^i\right) - n\sum_{i=0}^{k-1}(m-1)^i= \\
=  d\left(m\sum_{i=0}^{k-1}(m-1)^i + 1\right) -n\sum_{i=0}^{k-1}(m-1)^i.
\end{gather*}
\hfill $\qedsymbol$
\smallskip
\par\noindent
Next, we perform the final step of the computation, where our aim is to write the inequality that includes the term $X_z$. As in the previous cases, we divide into two cases.

\medskip
\par\noindent
\textbf{\underline{$r-1$ even:}} If $r-1 = 2\ell$, then we write the following inequality:
\begin{gather*}
    n \geq \dim{\left[\sum_{j=3}^m\left(\bigcap_{t=0}^{2\ell-1} X_{(j,t)}\right) + X_{(1,2\ell-2)} + X_{(1,2\ell)} + X_{(2,2\ell-1)}\right]}.
\end{gather*}
Due to the arguments used in the proof, the dimension can be decomposed into the following form:
\begin{gather*}
    n \geq (m-2)c_{2\ell-1} + 3d -\dim{\left[X_{(1,2\ell-2)}\cap X_{(1,2\ell)}\right]},
\end{gather*}
which, after rearranging, gives
\begin{gather*}
    \dim{\left[X_{(1,2\ell-2)}\cap X_{(1,2\ell)}\right]} \geq (m-2)c_{2\ell-1} + 3d -n.
\end{gather*}
Using the inequality calculated above, we derive an estimate for an intersection where every level appears except $t=2\ell-1$.
\begin{gather*}
    n \geq \dim{\left[\sum_{j=3}^m\left(\bigcap_{t=0}^{2\ell-3}X_{(j,t)}\right) + \bigcap_{t=0}^{2\ell-3}X_{(1,t)} + X_{(1,2\ell-2)}\cap X_{(1,2\ell)} + X_{(2,2\ell-3)}\right]}.
\end{gather*}
\smallskip
\par\noindent
It is worth noting that the appearance of the term $X_{(2,2\ell-3)}$ is the reason why the cases $r=2,3$ had to be treated separately.
\smallskip
\par\noindent
Observe that the term $X_{(2,2\ell-3)}$ can be factored out of the sum. Furthermore, for any $j=3,\dots,m$:
\begin{gather*}
    \sum_{\substack{s=3\\s\neq j}}^m\left(\bigcap_{t=0}^{2\ell-3}X_{(s,t)}\right) + \bigcap_{t=0}^{2\ell-3}X_{(1,t)} + X_{(1,2\ell-2)}\cap X_{(1,2\ell)} \subseteq \sum_{v \in N((j,2\ell-3))}X_v,
\end{gather*}
hence, the dimension can be decomposed as follows:
\begin{gather*}
    n \geq (m-1)c_{2\ell-3} + \dim{\left[X_{(1,2\ell-2)}\cap X_{(1,2\ell)}\right]} + d - \dim{\left[\bigcap_{\substack{t=0\\t\neq 2\ell-1}}^{2\ell}X_{(1,t)}\right]} \geq\\
    \geq (m-1)c_{2\ell-3} + (m-2)c_{2\ell-1} + 3d -n + d - \dim{\left[\bigcap_{\substack{t=0\\t\neq 2\ell-1}}^{2\ell}X_{(1,t)}\right]}.
\end{gather*}
From this, we obtain the following lower bound:
\[
\dim{\left[\bigcap_{\substack{t=0\\t\neq 2\ell-1}}^{2\ell}X_{(1,t)}\right]} \geq  (m-2)c_{2\ell-1} + (m-1)c_{2\ell-3} + 4d -2n.
\]
Finally, we consider an inequality that includes $X_z$:
\begin{gather*}
    n \geq \dim{\left[\sum_{j=2}^m\left(\bigcap_{\substack{t=0\\t \neq 2\ell-1}}^{2\ell}X_{(j,t)}\right) + X_{(1,2\ell-2)}\cap X_{(1,2\ell)} + X_z\right]}.
\end{gather*}
The term $X_z$ can be removed from the dimension count. Each remaining term has trivial intersection with the sum of the others. Thus,
\begin{gather*}
    n \geq \sum_{j=2}^m \dim{\left[\bigcap_{\substack{t=0\\ t \neq 2\ell-1}}^{2\ell}X_{(j,t)}\right]} +  \dim{\left(X_{(1,2\ell-2)}\cap X_{(1,2\ell)}\right)} + d \geq \\
    \geq (m-1)((m-2)c_{2\ell-1} + (m-1)c_{2\ell-3} + 4d-2n) +(m-2)c_{2\ell-1} + 4d-n. 
\end{gather*}
Rearranging, we obtain
\begin{gather}\label{ineq1}
    m(m-2)c_{r-2} + (m-1)^2c_{r-4} + 4dm - 2nm \leq 0.
\end{gather}
Before showing that this inequality leads to the correct final result, we first consider that we obtain the same inequality for the odd case.

\medskip
\par\noindent
\textbf{\underline{$r-1$ odd:}} Let $r-1 = 2\ell+1$ and proceed exactly as in the even case. First, we write:
\begin{gather*}
    n \geq \dim{\left[\sum_{j=3}^m\left(\bigcap_{t=0}^{2\ell} X_{(j,t)}\right) + X_{(1,2\ell-1)} + X_{(1,2\ell+1)} + X_{(2,2\ell)}\right]}.
\end{gather*}
This gives the estimate:
\begin{gather*}
    \dim{\left(X_{(1,2\ell-1)}\cap X_{(1,2\ell+1)}\right)} \geq (m-2)c_{2\ell} + 3d -n.
\end{gather*}
Next, we derive:
\begin{gather*}
    n \geq \dim{\left[\sum_{j=3}^m\left(\bigcap_{t=0}^{2\ell-2}X_{(j,t)}\right) + \bigcap_{t=0}^{2\ell-2}X_{(1,t)} + X_{(1,2\ell-1)}\cap X_{(1,2\ell+1)} + X_{(2,2\ell-2)}\right]},
\end{gather*}
from which we obtain:
\[
\dim{\left[\bigcap_{\substack{t=0\\t\neq 2\ell}}^{2\ell+1}X_{(1,t)}\right]} \geq  (m-2)c_{2\ell} + (m-1)c_{2\ell-2} + 4d -2n.
\]
The final inequality is the following:
\begin{gather*}
    n \geq \dim{\left[\sum_{j=2}^m\left(\bigcap_{\substack{t=0\\t \neq 2\ell}}^{2\ell+1}X_{(j,t)}\right) + X_{(1,2\ell-1)}\cap X_{(1,2\ell+1)} + X_z\right]}.
\end{gather*}
After expanding and rearranging, we obtain
\begin{gather*}
    m(m-2)c_{r-2} + (m-1)^2c_{r-4} + 4dm - 2nm \leq 0.
\end{gather*}
The final step in the proof is to verify the following identity:
\medskip
\par\noindent
\[
m(m-2)c_{r-2} + (m-1)^2c_{r-4} + 4dm - 2nm = c_r.
\]
Using the recurrence $c_{k+1}=(m-1)c_k+2d-n$ from Lemma~\ref{lem2}
we obtain by two substitutions that
\begin{gather*}
c_r=(m-1)c_{r-1}+2d-n
=\\
=(m-1)\bigl((m-1)c_{r-2}+2d-n\bigr)+2d-n
=(m-1)^2c_{r-2}+m(2d-n).
\end{gather*}
\par\noindent
Similarly,
\begin{gather*}
c_{r-2}=(m-1)c_{r-3}+2d-n
=\\
=(m-1)\bigl((m-1)c_{r-4}+2d-n\bigr)+2d-n
=(m-1)^2c_{r-4}+m(2d-n),
\end{gather*}
hence
\begin{gather*}
(m-1)^2c_{r-2}+m(2d-n) = m(m-2)c_{r-2} + c_{r-2} + m(2d-n) = \\
= m(m-2)c_{r-2}+(m-1)^2c_{r-4}+4dm-2nm.
\end{gather*}
\medskip
\par\noindent
Substituting this identity into inequality~(\ref{ineq1}) yields $c_r \le 0$, which is equivalent to the desired lower bound. This completes the proof of Theorem~\ref{Theorem2}.
\hfill $\qedsymbol$
\bigskip
\par\noindent
\section*{Acknowledgment}

The author would like to express his gratitude to his late advisor, G\'abor Simonyi, for his guidance and insight during the initial phase of this work.

\newpage
\section{Appendix}
\bigskip
\par\noindent
{\Large{\bf Appendix A}}
\bigskip
\par\noindent
{\em Proof of Lemma~\ref{intersection}.} The inclusion from right to left is trivial. We need to prove the other direction. Consider a vector
\begin{gather*}
     v \in\left[\sum_{i=1}^k\left(X_i\otimes \Gamma[a_i,b_i]\right)\right] \cap \left[\sum_{j=1}^{\ell}(Y_j\otimes\Gamma[c_j,d_j]) + Y \otimes \Gamma[c,d]\right].
\end{gather*}
\normalsize
For each $i$, let $x^{(1,i)},\dots,x^{(t_i,i)}$ be a basis of $X_i$. Similarly, for each $j$, let $y^{(1,j)},...,y^{(r_j,j)}$ be a basis of $Y_j$ and let $y^{(1)},...,y^{(r)}$ be a basis of $Y$. The vector $v$ can be expressed in two different ways corresponding to these bases:
\begin{gather*}
    v = \sum_{i=1}^k\sum_{h \in [a_i,b_i]}\sum_{p=1}^{t_i}\alpha_{p,h}x^{(p,i)}\otimes e^{(h)}.
\end{gather*}
and
\begin{gather*}
    v = \sum_{j=1}^{\ell}\sum_{s \in [c_j,d_j]}\sum_{q=1}^{r_j}\beta_{q,s}y^{(q,j)}\otimes e^{(s)} + \sum_{w \in [c,d]}\sum_{u=1}^r\gamma_{u,w} y^{(u)} \otimes e^{(w)},
\end{gather*}
where $e^{(h)},e^{(s)},e^{(w)}$'s are the canonical vectors of $\mathbb{F}^M$ and $\alpha_{p,h},\beta_{q,s},\gamma_{u,w} \in \mathbb{F}$. For simplicity we use the following notations:
\begin{gather*}
    \tilde{x}^{(h,i)}:= \sum_{p = 1}^{t_i}\alpha_{p,h}x^{(p,i)}, \hspace{0.5 cm}\tilde{y}^{(s,j)}:= \sum_{q=1}^{r_j}\beta_{q,s}y^{(q,j)}, \hspace{0.5 cm} \tilde{y}^{(w)}:=\sum_{u=1}^{r}\gamma_{u,w}y^{(u)}.
\end{gather*}
By the bilinearity of the tensor product, $v$ can be expressed as
\begin{gather*}
    v = \sum_{i=1}^k\sum_{h \in [a_i,b_i]}\tilde{x}^{(h,i)}\otimes e^{(h)};\\
    v = \sum_{j=1}^{\ell}\sum_{s \in [c_j,d_j]}\tilde{y}^{(s,j)}\otimes e^{(s)} + \sum_{w \in [c,d]}\tilde{y}^{(w)} \otimes e^{(w)}.
\end{gather*}
The intervals $[c_j,d_j]$ could intersect $[c,d]$, so we separate the second sum according to whether the indices lie in $[c,d]$ or outside.
\begin{gather*}
    \sum_{j=1}^{\ell}\sum_{s \in [c_j,d_j]\setminus [c,d{}]}\tilde{y}^{(s,j)}\otimes e^{(s)} + \sum_{j=1}^{\ell}\sum_{s \in [c,d]}\tilde{y}^{(s,j)}\otimes e^{(s)} + \sum_{w \in [c,d]}\tilde{y}^{(w)} \otimes e^{(w)}.
\end{gather*}
Define
\begin{gather*}\tilde{v}:=\sum_{i=1}^k\sum_{h \in [a_i,b_i]}\tilde{x}^{(h,i)}\otimes e^{(h)} -  \sum_{j=1}^{\ell}\sum_{s \in [c_j,d_j]\setminus [c,d{}]}\tilde{y}^{(s,j)}\otimes e^{(s)} = \\
    =\sum_{j=1}^{\ell}\sum_{s \in [c,d]}\tilde{y}^{(s,j)}\otimes e^{(s)} + \sum_{w \in [c,d]}\tilde{y}^{(w)} \otimes e^{(w)}.
\end{gather*}
Since  $\bigcup_{i=1}^k[a_i,b_i]$ and $\bigcup_{j=1}^{\ell}([c_j,d_j]\setminus [c,d])$ are disjoint from $[c,d]$, the intersection of the two spans is trivial, that is,
\begin{gather*}
    {\rm span}_{\mathbb{F}}\left\{e^{(h)}:h \in \bigcup_{i=1}^k[a_i,b_i]\cup\bigcup_{j=1}^{\ell}([c_j,d_j]\setminus [c,d])\right\} \cap {\rm span}_{\mathbb{F}}\left\{e^{(s)}:s \in [c,d]\right\} = \{0\}.
\end{gather*}
The above fact implies $\tilde{v} = 0$, so we have
\begin{gather*}
    v = \sum_{j=1}^{\ell}\sum_{s \in [c_j,d_j]\setminus [c,d]}\tilde{y}^{(s,j)}\otimes e^{(s)} \in \sum_{j=1}^{\ell}(Y_j\otimes\Gamma[c_j,d_j]).
\end{gather*}
\hfill $\qedsymbol$
\par\noindent
{\Large{\bf Appendix B1}}
\bigskip
\par\noindent
In this appendix, we verify that for every $u \in V(M_r(G))$ the property $\tilde{X}_u \cap \left(\sum_{w \in N_{M_r(G)}(u)}\tilde{X}_w\right) = \{0\}$ holds. As shown in the proof of Theorem~\ref{Theorem1}, the property holds for vertices of the form $(v,0)$.
\medskip
\par\noindent
\textbf{\underline{$r = 2$:}} Let $v \in V(G)$ be arbitrary.
\footnotesize
\begin{gather*}
    \tilde{X}_{(v,1)} \cap \left(\sum_{u \in N_{M_2(G)}((v,1))}\tilde{X}_u\right) =  \tilde{X}_{(v,1)} \cap \left(\sum_{q \in N_G(v)}\tilde{X}_{(q,0)} + \tilde{X}_z\right) =\\ 
    = \left((X_v \otimes \Gamma[d^2+1,nd]+E\otimes\Gamma[nd+1,nd + d^3]\right) \cap \left[\sum_{q \in N_G(v)}\left( X_q \otimes \Gamma[1,nd]\right) + \left(\sum_{u \in V(G)}X_u\right)\otimes \Gamma[1,d^2]\right].
\end{gather*}
\normalsize
Using Lemma~\ref{intersection}, we may disregard the terms $\left(\sum_{u \in V(G)}X_u\right)\otimes \Gamma[d^2+1,nd]$ and $E\otimes\Gamma[nd+1,nd + d^3]$. Thus,
\smallskip
\par\noindent
\small
\begin{gather*}
    (X_v \otimes \Gamma[d^2 +1,nd]) \cap \left[\sum_{q \in N_G(v)} (X_q \otimes \Gamma[1,nd])\right] \subseteq (X_v \otimes \Gamma[1,nd]) \cap \left[\left(\sum_{q \in N_G(v)} X_q\right) \otimes \Gamma[1,nd]\right].
\end{gather*}
\normalsize
Using Theorem~\ref{tensor}, the intersection is trivial.
\medskip
\par\noindent
Now we check the adjacency property of $\tilde{X}_z$.
\begin{gather*}
    \tilde{X}_z \cap \left(\sum_{u \in N_{M_2(G)}(z)}\tilde{X}_u\right)=\\ = \left[\left(\sum_{v \in V(G)}X_v\right) \otimes \Gamma[1,d^2]\right] \cap \left[\sum_{v \in V(G)}\left(X_v \otimes \Gamma[d^2 +1,nd]\right) + E\otimes\Gamma[nd+1,nd+d^3]\right] = \{0\}.
\end{gather*}
Since $[1,d^2]\cap[d^2+1,nd+d^3]=\emptyset$, Lemma~\ref{intersection} implies that the intersection is trivial.
\medskip
\par\noindent
\textbf{\underline{Case of $r = 2t+1$:}} We first verify the adjacency condition for $(v,1)$.
\small
\begin{gather*}
    \tilde{X}_{(v,1)} \cap \left(\sum_{u \in N_{M_r(G)}((v,1))}\tilde{X}_u\right) = \tilde{X}_{(v,1)} \cap \left(\sum_{q \in N_G(v)}\tilde{X}_{(q,0)} + \sum_{q \in N_G(v)}\tilde{X}_{(q,2)}\right) = \\
    = \left(X_v \otimes \Gamma[d^{2r-2}+1,a_{r-1}] + E\otimes\Gamma\left[a_{r-1} + 1, a_{r-1} + d^{2r-1}\right]\right) \cap \\
    \cap \left(\sum_{q \in N_G(v)}\left(X_q \otimes \Gamma[1,a_{r-1}]\right) + \left(\sum_{w\in V(G)}X_w\right) \otimes \Gamma\left[1,d^{2r-2}\right]+\sum_{q \in N_G(v)}\left(X_q \otimes\Gamma[a_1 +1,a_{r-1}]\right)\right).
\end{gather*}
\normalsize
We can disregard the terms $E \otimes\Gamma\left[a_{r-1} + 1, a_{r-1} + d^{2r-1}\right],\left(\sum_{w\in V(G)}X_w\right) \otimes \Gamma\left[1,d^{2r-2}\right]$ and we have
\smallskip
\par\noindent
\begin{gather*}
    \left(X_v \otimes \Gamma[d^{2r-2}+1,a_{r-1}] \right) \cap \\
    \cap \left(\sum_{q \in N_G(v)}\left(X_q \otimes \Gamma[1,a_{r-1}]\right) + \sum_{q \in N_G(v)}\left(X_q \otimes \Gamma[a_1 +1,a_{r-1}]\right)\right)\subseteq\\
    \subseteq
    \left(X_v \otimes \Gamma[1,a_{r-1}]\right) \cap \left[\left(\sum_{q \in N_G(v)}X_q\right) \otimes \Gamma[1,a_{r-1}]\right] = \{0\}.
\end{gather*}
Let $k=2\ell < r-1$.
\smallskip
\par\noindent
\scriptsize
\begin{gather*}
    \tilde{X}_{(v,k)}\cap\left(\sum_{u \in N_{M_r(G)}((v,k))}\tilde{X}_u\right) = \tilde{X}_{(v,k)} \cap\left(\sum_{q \in N_G(v)}\tilde{X}_{(q,k-1)} + \sum_{q \in N_G(v)}\tilde{X}_{(q,k+1)}\right) =\\
    = \left[\left(\sum_{w \in V(G)}X_w\right) \otimes \left(\sum_{i=0}^{\ell-1} \Gamma[a_{2i-1}+1,a_{2i}]\right) + X_v \otimes \Gamma[a_{k-1}+1,a_{r-1}]\right]\cap\\
    \cap \left[\left(\sum_{w \in V(G)}X_w\right) \otimes\left(\sum_{i=0}^{\ell-1} \Gamma[a_{2i}+1,a_{2i+1}]\right) + \left(\sum_{q \in N_G(v)}X_q\right)\otimes \Gamma[a_{k-2}+1,a_{r-1}] + E\otimes\Gamma\left[a_{r-1}+1,a_{r-1}+d^{2r-1}\right]\right].
\end{gather*}
\normalsize
We may disregard the term $E\otimes \Gamma[ar-1+1,ar-1+d^{2r-1}]$. The remaining intersection is trivial by Theorem~\ref{tensor}.
\medskip
\par\noindent
Let $k=2\ell+1 < r-1$. 
\smallskip
\par\noindent
\footnotesize
\begin{gather*}
\tilde{X}_{(v,k)}\cap\left(\sum_{u \in N_{M_r(G)}((v,k))}\tilde{X}_u\right) = \tilde{X}_{(v,k)} \cap\left(\sum_{q \in N_G(v)}\tilde{X}_{(q,k-1)} + \sum_{q \in N_G(v)}\tilde{X}_{(q,k+1)}\right) = \\
= \left[\left(\sum_{w \in V(G)}X_w\right) \otimes\left(\sum_{i=0}^{\ell-1} \Gamma[a_{2i}+1,a_{2i+1}]\right) + X_v\otimes\Gamma[a_{k-1}+1,a_{r-1}] + E\otimes\Gamma\left[a_{r-1}+1,a_{r-1}+d^{2r-1}\right]\right]\cap\\
\cap \left[\left(\sum_{w \in V(G)}X_w\right)\otimes\left(\sum_{i=0}^{\ell}\Gamma[a_{2i-1}+1,a_{2i}]\right)+\left(\sum_{ q \in N_G(v)}X_q\right)\otimes\Gamma[a_{k-2}+1,a_{r-1}]\right].
\end{gather*}
\normalsize
We can disregard the term $E\otimes\Gamma[a_{r-1}+1,a_{r-1}+d^{2r-1}]$ and use Theorem~\ref{tensor}.
\medskip
\par\noindent
Let $k=2t = r-1$.
\scriptsize
\smallskip
\par\noindent
\begin{gather*}
    \tilde{X}_{(v,k)}\cap\left(\sum_{u \in N_{M_r(G)}((v,k))}\tilde{X}_u\right) = \tilde{X}_{(v,k)} \cap\left(\sum_{q \in N_G(v)}\tilde{X}_{(q,k-1)} + \tilde{X}_z\right) \subseteq \\
    \subseteq \left[ \left(\sum_{w \in V(G)}X_w\right) \otimes \left(\sum_{i=0}^{t-1} \Gamma[a_{2i-1}+1,a_{2i}]\right)+X_v \otimes \Gamma[a_{r-2}+1,a_{r-1}]\right]\cap\\
    \cap \left[\left(\sum_{w \in V(G)}X_w \right)\otimes\left(\sum_{i=0}^{t-1} \Gamma[a_{2i}+1,a_{2i+1}] \right) + \left(\sum_{q \in N_G(v)}X_q \right)\otimes\Gamma[a_{r-3}+1,a_{r-1}] +E\otimes\Gamma\left[a_{r-1}+1,a_{r-1}+d^{2r-1}\right]\right].
\end{gather*}
\normalsize
\par\noindent
We can disregard the term $E\otimes\Gamma\left[a_{r-1}+1,a_{r-1}+d^{2r-1}\right]$.
\par\noindent
Furthermore $\left(\sum_{i=0}^{t-1} \Gamma[a_{2i-1}+1,a_{2i}]\right) \cap \left(\sum_{i=0}^{t-1} \Gamma[a_{2i}+1,a_{2i+1}] \right) = \{0\}$. Hence, any vector in the intersection must lie in $\Gamma[ar-2+1,ar-1]$. Theorem~\ref{tensor} then implies that the intersection is trivial.
\medskip
\par\noindent
Now we show that for $\tilde{X}_z$ the adjacency property holds.
\begin{gather*}
    \tilde{X}_z \cap \left(\sum_{v \in V(G)}\tilde{X}_{(v,r-1)}\right) = \\
    = \left[\left(\sum_{v \in V(G)}X_v\right)\otimes\left(\sum_{i=0}^{t-1} \Gamma[a_{2i}+1,a_{2i+1}]\right) + E \otimes\Gamma\left[a_{r-1}+1,a_{r-1}+d^{2r-1}\right]\right] \cap \\
    \cap  \left[\left(\sum_{v \in V(G)}X_v\right) \otimes \left(\sum_{i=0}^{t-1} \Gamma[a_{2i-1}+1,a_{2i}]+\Gamma[a_{r-2}+1,a_{r-1}]\right)\right].
\end{gather*}
Apply Theorem~\ref{tensor}.
\medskip
\par\noindent
\textbf{\underline{Case of $r=2t$:}} In this case, for every $k=0,1,..,r-2$ the adjacency property of $(v,k)$ is similar. Assume that $k=2t-1 =r-1$.
\smallskip
\par\noindent
\scriptsize
\begin{gather*}
    \tilde{X}_{(v,k)} \cap \left(\sum_{q \in N_G(v)}\tilde{X}_{(q,k-1)} + \tilde{X}_z\right) =\\
    =\left[\left(\sum_{w \in V(G)}X_w\right) \otimes\left(\sum_{i=0}^{t-2} \Gamma[a_{2i}+1,a_{2i+1}]\right) + X_v \otimes \Gamma[a_{r-2}+1,a_{r-1}] + E\otimes\Gamma\left[a_{r-1}+1,a_{r-1}+d^{2r-1}\right]\right]\cap \\
    \cap \left[\left(\sum_{w \in V(G)}X_w\right) \otimes \left(\sum_{i=0}^{t-1} \Gamma[a_{2i-1}+1,a_{2i}]\right)+\left(\sum_{q \in N_G(v)}X_q\right)\otimes \Gamma[a_{r-3}+1,a_{r-1}]\right]
\end{gather*}
\normalsize
As in the case $r=2t+1$, the intersection is trivial.
\medskip
\par\noindent
Finally, we consider the subspace $\tilde{X}_z$.
\smallskip
\par\noindent
\footnotesize
\begin{gather*}
     \tilde{X}_z \cap \left(\sum_{v \in V(G)}\tilde{X}_{(v,r-1)}\right)    = \left[\left(\sum_{w \in V(G)}X_w\right)\otimes\left(\sum_{i=0}^{t-1} \Gamma[a_{2i-1}+1,a_{2i}]\right)\right] \cap \\
    \cap \left[\left(\sum_{v \in V(G)}X_v\right) \otimes\left(\sum_{i=0}^{t-1} \Gamma[a_{2i}+1,a_{2i+1}] + \Gamma[a_{r-2}+1,a_{r-1}]\right) + E\otimes\Gamma\left[a_{r-1}+1,a_{r-1}+d^{2r-1}\right]\right].
\end{gather*}
\normalsize
This intersection is trivial since 
\smallskip
\par\noindent
\scriptsize
\[
\left(\sum_{i=0}^{t-1} \Gamma[a_{2i-1}+1,a_{2i}]\right) \cap \left(\sum_{i=0}^{t-1} \Gamma[a_{2i}+1,a_{2i+1}] + \Gamma[a_{r-2}+1,a_{r-1}] + E\otimes\Gamma\left[a_{r-1}+1,a_{r-1}+d^{2r-1}\right]\right) = \{0\}.
\]
\normalsize
This completes the verification of the adjacency condition for the constructed dual representation.
\bigskip
\par\noindent
{\Large{\bf Appendix B2}}
\bigskip
\par\noindent 
In this appendix, we verify that all subspaces assigned $\tilde{X}_{(v, k)}$ and $\tilde{X}_z$ have the required dimension $D$. In the proof, we use the identity $\dim{(X \otimes Y)} = \dim{X}\dim{Y}$ from Theorem~\ref{tensor}.
\smallskip
\par\noindent
\begin{gather*}
    \dim{\tilde{X}_{(v,0)}} = \dim{\left(X_v \otimes\Gamma[1,a_{r-1}]\right)} = \dim{X_v}\dim{\Gamma[1,a_{r-1}]}=d\sum_{i=0}^{r-1}d^{2r-2-i}(n-d)^i = D.
\end{gather*}
\begin{gather*}
    \dim{\tilde{X}_{(v,1)} = \dim{\left( X_v \otimes\Gamma[d^{2r-2}+1,a_{r-1}] + E\otimes\Gamma\left[a_{r-1}+1,a_{r-1}+d^{2r-1}\right]\right)}} = \\
    =d\left(a_{r-1} -d^{2r-2}\right) + d^{2r-1} = D.
\end{gather*}
Let $k=2\ell$.
\smallskip
\par\noindent
\footnotesize
\begin{gather*}
    \dim\tilde{X}_{(v,k)}
    =\dim{\left[\left(\sum_{w \in V(G)}X_w\right) \otimes \left(\sum_{i=0}^{\ell-1} \Gamma[a_{2i-1}+1,a_{2i}]\right)+X_v \otimes\Gamma[a_{k-1}+1,a_{r-1}]\right]}
    =\\
    =n\sum_{i=0}^{\ell-1}(a_{2i}-a_{2i-1}) + d(a_{r-1}-a_{2\ell-1}) = n\sum_{j=0}^{\ell-1}d^{2r-2-2j}(n-d)^{2j} + d\sum_{t=k}^{r-1}d^{2r-2-t}(n-d)^{t} = \\
    =nd^{2r-2}\sum_{j=0}^{\ell-1}\left(\frac{n-d}{d}\right)^{2j}+d^{2r-1}\sum_{t=k}^{r-1}\left(\frac{n-d}{d}\right)^{t} 
    =nd^{2r-2}\frac{\left(\frac{n-d}{d}\right)^{2\ell}-1}{\left(\frac{n-d}{d}\right)^2-1}+\\
    +d^{2r-1}\left(\frac{n-d}{d}\right)^{2\ell}\frac{\left(\frac{n-d}{d}\right)^{r-2\ell}-1}{\left(\frac{n-d}{d}\right)-1} 
    =d^{2r-2}\left(n\frac{\left(\frac{n-d}{d}\right)^{2\ell}-1}{\left(\frac{n-d}{d}\right)^2-1} + d\frac{\left(\frac{n-d}{d}\right)^r-\left(\frac{n-d}{d}\right)^{2\ell}}{\left(\frac{n-d}{d}\right)-1}\right)=\\
    =d^{2r-2}\left(\frac{n\left(\frac{n-d}{d}\right)^{2\ell}-n + d\left(\frac{n-d}{d}\right)^{r+1} - d\left(\frac{n-d}{d}\right)^{2\ell+1}+d\left(\frac{n-d}{d}\right)^r -d \left(\frac{n-d}{d}\right)^{2\ell}}{\left(\frac{n-d}{d}\right)^2-1}\right)=\\
    =d^{2r-2}n\frac{\left(\frac{n-d}{d}\right)^r-1}{\left(\frac{n-d}{d}\right)^2-1} =
    d^{2r-1}\frac{\left(\frac{n-d}{d}\right)^{r}-1}{\left(\frac{n-d}{d}\right)-1}= d^{2r-1}\sum_{i=0}^{r-1}\left(\frac{n-d}{d}\right)^{i} = D.
\end{gather*}
\normalsize
Let $k=2\ell+1$.
\smallskip
\par\noindent
\scriptsize
\begin{gather*}
    \dim{\tilde{X}_{(v,k)}} = \\
    =\dim{\left[\left(\sum_{w \in V(G)}X_w\right) \otimes\left(\sum_{i=0}^{\ell-1} \Gamma[a_{2i}+1,a_{2i+1}]\right) + X_v \otimes\Gamma[a_{k-1}+1,a_{r-1}] + E\otimes\Gamma\left[a_{r-1}+1,a_{r-1}+d^{2r-1}\right]\right]} =\\
    = n\left(\sum_{i=0}^{\ell-1}(a_{2i+1}-a_{2i})\right) + d(a_{r-1} -a_{2\ell}) = n\sum_{j=0}^{\ell-1}d^{2r-2-(2j+1)}(n-d)^{2j+1} + d\sum_{t=k}^{r-1}d^{2r-2-t}(n-d)^t +d^{2r-1}= \\
    = d^{2r-2}\left[n\left(\frac{n-d}{d}\right)\sum_{j=0}^{\ell-1}\left(\frac{n-d}{d}\right)^{2j} +d\sum_{t=k}^{r-1}\left(\frac{n-d}{d}\right)^t + d\right] = \\
    = d^{2r-2}\left[n\left(\frac{n-d}{d}\right)\frac{\left(\frac{n-d}{d}\right)^{2\ell}-1}{\left(\frac{n-d}{d}\right)^2-1} +\left(\frac{n-d}{d}\right)^{2\ell+1}\frac{\left(\frac{n-d}{d}\right)^{r-2\ell-1}-1}{\left(\frac{n-d}{d}\right)-1} + d\right] =\\
    = d^{2r-2}\left[\frac{n\left(\frac{n-d}{d}\right)^{2\ell+1}-n\left(\frac{n-d}{d}\right) + d\left(\frac{n-d}{d}\right)^{r+1}-d\left(\frac{n-d}{d}\right)^{2\ell+2}+d\left(\frac{n-d}{d}\right)^r -d\left(\frac{n-d}{d}\right)^{2\ell+1}}{\left(\frac{n-d}{d}\right)^2-1}\right]= \\
    =d^{2r-2}n\frac{\left(\frac{n-d}{d}\right)^r-1}{\left(\frac{n-d}{d}\right)^2-1} =
    d^{2r-1}\frac{\left(\frac{n-d}{d}\right)^{r}-1}{\left(\frac{n-d}{d}\right)-1}= d^{2r-1}\sum_{i=0}^{r-1}\left(\frac{n-d}{d}\right)^{i} = D.
\end{gather*}
\normalsize
Assume that $r = 2t$ and consider the dimension of $\tilde{X}_z$.
\smallskip
\par\noindent
\begin{gather*}
    \dim{\tilde{X}_z} = \dim{\left[\left(\sum_{v \in V(G)}X_v\right)\otimes\left(\sum_{i=0}^{t-1} \Gamma[a_{2i-1}+1,a_{2i}]\right)\right]} = n\sum_{i=0}^{t-1}(a_{2i}-a_{2i-1}) = \\
    = n\sum_{j=0}^{t-1}d^{2r-2-2j}(n-d)^{2j} = n d^{2r-2}\sum_{j=0}^{t-1}\left(\frac{n-d}{d}\right)^{2j} = nd^{2r-2}\frac{\left(\frac{n-d}{d}\right)^{r}-1}{\left(\frac{n-d}{d}\right)^2-1} = D.
\end{gather*}
Assume that $r=2t+1$.
\smallskip
\par\noindent
\small
\begin{gather*}
    \dim{\tilde{X}_z} = \dim{\left[\left(\sum_{v \in V(G)}X_v\right)\otimes\left(\sum_{i=0}^{t-1} \Gamma[a_{2i}+1,a_{2i+1}]\right) + E\otimes\Gamma\left[a_{r-1}+1,a_{r-1}+d^{2r-1}\right]\right]} = \\
    =n\sum_{i=0}^{t-1}(a_{2i+1}-a_{2i}) + d^{2r-1} = n\sum_{j=0}^{t-1}d^{2r-2-(2j+1)}(n-d)^{2j+1} + d^{2r-1} = \\
    = nd^{2r-2}\left(\frac{n-d}{d}\right)\sum_{j=0}^{t-1}\left(\frac{n-d}{d}\right)^{2j} + d^{2r-1} = nd^{2r-2}\left(\frac{n-d}{d}\right)\frac{\left(\frac{n-d}{d}\right)^{r-1}-1}{\left(\frac{n-d}{d}\right)^2-1} + d^{2r-1}=D.
\end{gather*}
\normalsize
This completes the verification that all subspaces in the construction have dimension \( D \), as required for a valid dual representation.

\begin{thebibliography}{14}

\bibitem{Anna}
Anna Blasiak, A graph-theoretic approach to network coding. Ph.D. thesis, {\em Cornell University}, (2013).

\bibitem{Boris}
Boris Bukh, Cristopher Cox,
On a Fractional Version of Haemers’ Bound, {\em IEEE Trans. Inform. Theory}, {\bf 65} (2019), 3340 -- 3348. 

\bibitem{Csonka}
Bence Csonka, G\'abor Simonyi,
Shannon Capacity, Lov\'asz Theta Number and the Mycielski Construction, {\em IEEE Trans. Inform. Theory}, {\bf 70} (2024), 7632 -- 7646. 

\bibitem{Imre11} Imre Csisz\'ar and J\'anos Körner, {\em Information Theory: Coding Theorems for Discrete Memoryless Systems. Second edition.}, Cambridge University Press, (2011).

\bibitem{Haemers}
Willem Haemers,
On some problems of Lov\'asz concerning the Shannon capacity of a graph,
{\em IEEE Trans. Inform. Theory}, {\bf 25} (1979),
231--232.

\bibitem{Haemers2}
Willem Haemers,
An upper bound for the Shannon capacity of a graph, {\em Algebraic Methods in Graph Theory}, Vol. {\bf 25} (1978), 267–-272.

\bibitem{Hoffman}
Kenneth Hoffman, Ray Kunze, Linear Algebra. 2nd edition, {\em Prentice Hall}, (1971).

\bibitem{LPU} Michael Larsen, James Propp, Daniel Ullman, The fractional chromatic number of Mycielski graphs, {\em J. Graph Theory} {\bf 19} (1995), 411--416.

\bibitem{LL79} L\'aszl\'o Lov\'asz,
On the Shannon capacity of a graph,
{\em IEEE Trans. Inform. Theory }, {\bf 25} (1979),
1--7.

\bibitem{Myc} Jan Mycielski, Sur le coloriage des graphs, {\em Colloq.\ Math.},
{\bf 3} (1955), 161--162.

\bibitem{Ryan} Raymond A. Ryan, Introduction to Tensor Products of Banach Spaces {\em Springer Monographs in Mathematics},
(2002).

\bibitem{Sha56}
Claude E. Shannon,
The zero-ERROR capacity of a noisy channel,
{\em IRE Trans. Inform. Theory}, {\bf 2} (1956), 8--19.

\bibitem{CTcone} Claude Tardif, Fractional chromatic numbers of cones over graphs,
{\em J. Graph Theory} {\bf 38} (2001), 87--94.

\bibitem{Zuiddam}
Jeroen Zuiddam, The Asymptotic Spectrum of Graphs and the Shannon Capacity, {\em Combinatorica}, {\bf 39} (2019), 1173–-1184.
\end{thebibliography}
\end{document}